\documentclass{ifacconf}

\usepackage{graphicx}      
\usepackage{natbib}        
\usepackage{cite}
\usepackage{amsmath,amssymb,amsfonts}
\usepackage{algorithmic}
\usepackage{graphicx}
\usepackage{textcomp}
\usepackage{xcolor}

\usepackage{latexsym, amssymb, amsmath}
\usepackage{enumerate}
\usepackage{flushend} 
\usepackage{mathrsfs} 
\usepackage{stfloats}
\def\custombibliography#1{
 \normalsize
\section*{\centering References}
 \list
 {[\arabic{enumi}]}{\settowidth\labelwidth{[#1]}\leftmargin\labelwidth
 \setlength{\itemsep}{.1em}
 \advance\leftmargin\labelsep
 \usecounter{enumi}}
 \def\newblock{\hskip .11em plus .33em minus -.07em}
 \sloppy
 \sfcode`\.=1000\relax}

\def\L2{{\cal L}_2}

\newcommand\bull{\vrule height .9ex width .8ex depth -.1ex } 

\newcommand\re{\rm I\! R}
\newcommand\cdcout[1]{} 

\newcommand{\rv}[1]{\boldsymbol{#1}} 
\newcommand{\RomanNumber}[1]{\uppercase\expandafter{\romannumeral #1}}
\newcommand{\romannumber}[1]{\lowercase\expandafter{\romannumeral #1}}

\DeclareMathAlphabet{\mathpzc}{OT1}{pzc}{m}{it}

\def\1{\rv 1} 

\usepackage{flushend} 
\usepackage{stfloats}
\usepackage{color} 
\usepackage{latexsym,amssymb,amsmath} 
\usepackage{mathrsfs,mathtools} 
\usepackage{multicol,enumerate}
\usepackage{subcaption}

\def\abs#1{\lvert #1 \rvert}

\def\allcommutingseries{\mbox{$\re\, [[\tilde{X}]]$}}

\def\allcommutingseriesXell{\mbox{$\re^{\ell}\, [[X]]$}} 

\def\allcommutingseriesX#1{\mbox{$\re^{#1}\, [[\tilde{X}]]$}}
\def\allcommutingpiseries#1{\mbox{$\re^{#1}_{pi}\, [[\tilde{X}]]$}}
\def\allcommutingseriesLCX#1{\mbox{$\re_{LC}^{#1}\, [[\tilde{X}]]$}}

\def\allpolyx0degn{\mbox{$P_n$}}

\def\allwords{\mbox{$X^{\ast}$}}
\def\allcommutingwords{\mbox{$\tilde{X}^{\ast}$}}
\def\allseries{\mbox{$\re\langle\langle X \rangle\rangle$}}
\def\allseries#1{\mbox{$\re^{#1}\langle\langle X \rangle\rangle$}}
\def\allseriesXpri#1{\mbox{$\re^{#1}\langle\langle X^{\prime} \rangle\rangle$}}

\def\allseriesell{\mbox{$\re^{\ell} \langle\langle X \rangle\rangle$}}

\def\allproperseriesell{\mbox{$\re_{p}^{\ell}\, \langle\langle X \rangle\rangle$}}
\def\allproperseries#1{\mbox{$\re_{p}^{#1}\, \langle\langle X \rangle\rangle$}}
\def\allnproperseries#1{\mbox{$\re_{np}^{#1}\, \langle\langle X \rangle\rangle$}}
\def\allpiseries#1{\mbox{$\re_{pi}^{#1}\, \langle\langle X \rangle\rangle$}}
\def\allpiseriesXpri#1{\mbox{$\re_{pi}^{#1}\, \langle\langle X^{\prime} \rangle\rangle$}}
\def\de#1{\mbox{${}^{\footnotesize{{\delta}}}#1$}}

\def\deltaFliess{\mbox{${}^{\delta}\mathscr{F}$}}

\def\allseriesellLC{\mbox{$\re^{\ell}_{LC}\langle\langle X \rangle\rangle$}}

\def\allseriesX1{\mbox{$\re [[ X_1 ]]$}}

\def\bull{\rule{0.08in}{0.08in}} 

\newcommand{\comment}[1]{} 

\def\doubleone{{\rm\, l\!l}}

\def\Endallseries{{\rm End}\left(\allseries{}\right)}

\def\Homstat#1#2{{\rm Hom}_{\tiny{\rm static}} \left(\re^{#1},\re^{#2}\right)}
\def\eqref#1{(\ref{#1})} 

\def\mbf#1{\hbox{\mathversion{bold}$#1$}} 
\def\mixcomp{\:\hat{\circ}\,} 
\def\mixprocomp{\:\check{\circ}\,} 
\def\nat{{\mathbb N}} 
\def\norm#1{\Vert#1\Vert}

\def\openbull{\framebox[0.08in][c]{$\;$}} 
\def\ord{{\rm ord}}

\def\re{{\mathbb R}} 

\def\shuffle{{\scriptscriptstyle \;\sqcup \hspace*{-0.05cm}\sqcup\;}}

\def\supp{{\rm supp}}
\def\orde{{\rm ord}}

\def\vec{{\rm vec}}

\def\begals{\[\begin{aligned}}
\def\endals{\end{aligned}\]}
\def\begal{\begin{align*}}
\def\endal{\end{align*}}
\def\begce{\begin{center}}
\def\endce{\end{center}}
\def\begar{\begin{array}}
\def\endar{\end{array}}
\def\begeq{\begin{equation}}
\def\endeq{\end{equation}}
\def\begdi{\begin{displaymath}}
\def\enddi{\end{displaymath}}
\def\begdis{\begin{eqnarray*}}
\def\enddis{\end{eqnarray*}}
\def\begeqa{\begin{eqnarray}}
\def\endeqa{\end{eqnarray}}
\def\begdes{\begin{description}}
\def\enddes{\end{description}}
\def\begit{\begin{itemize}}
\def\endit{\end{itemize}}
\def\begen{\begin{enumerate}}
\def\enden{\end{enumerate}}
\def\beglar{\left[\begin{array}}
\def\endrar{\end{array}\right]}
\def\begle{\begin{lemma}}
\def\endle{\end{lemma}}
\def\begde{\begin{definition}}
\def\endde{\end{definition}}
\def\begth{\begin{theorem}}
\def\endth{\end{theorem}}
\def\begco{\begin{corollary}}
\def\endco{\end{corollary}}
\def\begprop{\begin{proposition}}
\def\endprop{\end{proposition}}
\def\begex{\begin{example}}
\def\endex{\hfill\openbull \end{example}}
\def\begexer{\begin{exercise}}
\def\endexer{\end{exercise}}
\def\begalg{\begin{algo}}
\def\endalg{\end{algo}}

\def\begres{\noindent{\bf Remarks}:\begin{enumerate}}
\def\endres{\end{enumerate} \par}
\def\begpr{\noindent{\em Proof:}$\;\;$}
\def\endpr{\hfill\bull}
\def\begtab{\begin{tabular}}
\def\endtab{\end{tabular}}
\def\rref#1{(\ref{#1})}
\allowdisplaybreaks[4]

\newtheorem{lemma}{Lemma}[section]
\newtheorem{definition}{Definition}[section]
\newtheorem{theorem}{Theorem}[section]
\newtheorem{corollary}{Corollary}[section]
\newtheorem{example}{Example}[section]

\begin{document}
	
\begin{frontmatter}

\title{A Formal Power Series Approach to Multiplicative Dynamic and Static Output Feedback} 

\thanks[footnoteinfo]{Corresponding Author}
\thanks{accepted in $25^{th}$ Int. Symposium on Mathematical Theory of Networks and Systems, Bayreuth, Germany, 2022.}
\author[First]{Venkatesh~G.~S.} 

\address[First]{Old Dominion University, 
   Norfolk, VA 23529 USA (e-mail: sgugg001@odu.edu).}

\begin{abstract}                
The goal of the paper is two-fold. The first of which is to derive an explicit formula to compute the generating series of a closed-loop system when a plant, given in a Chen-Fliess series description is in multiplicative output feedback connection with another system given in Chen-Fliess series description. Further, the objective extends in showing that the multiplicative dynamic output feedback connection has a natural interpretation as a transformation group acting on the plant. The second of the two-part goal of this paper is same as the first part albeit when the Chen-Fliess series in the feedback is replaced by a memoryless map. The paper provides an explicit formula to compute the generating series of a closed-loop system in multiplicative static output feedback connection and shows that the static feedback has a natural interpretation as a transformation group acting on the plant.       	
\end{abstract}


\end{frontmatter}

\section{Introduction}
The objective of the document is two fold and works with the Chen-Fliess functional series~\citep{Fliess_81}. There is no need that these input-output systems have a state space realization and thus, the results presented here are independent of any state space embedding when a realization is possible~\citep{Fliess_realizn_83}. Firstly, let $F_c$ and $F_d$ be two nonlinear input-output systems represented by Chen-Fliess series. It was shown in~\citet{Gray-Li_05} that the {\em additive feedback} interconnection of two such systems result in a Chen-Fliess series description for the closed-loop system. An efficient computation of the generating series for closed-loop system is facilitated through a combinatorial Hopf algebra~\citep{Gray-etal_SCL14,Duffaut-Espinosa-etal_JA16}. The convergence of the closed-loop system was characterized in~\citep{Thitsa-Gray_SIAM12}. The feedback product formula and its computation were used to solve system inversion problems~\citep{Gray-etal_AUTO14} and trajectory generation problems~\citep{Duffaut-Espinosa-Gray_ICSTCC17}. However, when the nature of interconnection becomes {\em multiplicative feedback}, the similar set of questions persist in general. It is known that, in single-input single-output setting, the closed-loop system in the affine feedback case (of which multiplicative feedback is a special case) has a Chen-Fliess series description and the computation of feedback formula is facilitated through a combinatorial Hopf algebra~\citep{Gray-KEF_SIAM2017}. The present document, in one part, shows that even in multi-input multi-output setting the closed-loop system under multiplicative feedback has a Chen-Fliess series representation and provides an explicit expression of the closed-loop generating series and will be called as {\em multiplicative dynamic feedback product} . It will be shown that this feedback product has a natural interpretation as a transformation group acting on the plant. The algorithmic framework for the computation of the multiplicative dynamic feedback product formula for a general multi-input multi-output case and characterization of convergence for the closed-loop system is deferred to future work. Hence, the document is void of a computational example.

Secondly, let $F_d$ in the feedback path be replaced by a memoryless function $f_d$ which is coined as {\em additive static feedback} connection~\citep{Isidori_95}. The closed-loop system for the additive static feedback interconnection is known to have a Chen-Fliess series representation and an explicit expression for the closed-loop generating series, called {\em Wiener-Fliess feedback product}~\citep{Venkat-Gray_2021}. An algorithmic framework for computing the feedback product through the interaction of two Hopf algebras, Hopf algebra of the shuffle group and Hopf algebra of the dynamic feedback group, is presented in \citet{Venkat-Gray_2021}. The convergence of the closed-loop system was characterized in ~\citet{GS_thesis}. However the questions remain open when the nature of static feedback becomes {\em multiplicative}. Hence, the second of the two-part goal of this paper is to show that the closed-loop system in {multiplicative static feedback connection} has a Chen-Fliess series representation and an explicit expression for the closed-loop generating series, called {\em multiplicative static feedback product}, is provided. Further, the feedback product is shown as a transformation group acting on the plant. As in the case of multiplicative dynamic feedback product, the algorithmic framework for the computation of the multiplicative static feedback product and characterization of convergence for the closed-loop system is deferred to future work.            
 
 The paper is organized as follows. The next section provides a summary of the concepts related to Chen-Fliess series and their interconnections. The section also build the pivotal {\em multiplicative dynamic output feedback group} and also provides a brief discussion on formal static maps and Wiener-Fliess composition. Section~\ref{sec:mult_dyn_feedb} is where the multiplicative dynamic feedback connection is analyzed an Section~\ref{sec:mult_stat_feedb} is where the results of the multiplicative static feedback connection is detailed. The conclusions of the paper and directions for future work is given in the last section.
 
\section{Preliminaries}
A finite nonempty set of noncommuting symbols $X=\{ x_0,x_1,\ldots,x_m\}$ is called an {\em alphabet}. Each element of $X$ is called a {\em letter}, and any finite sequence of letters from $X$, $\eta=x_{i_1}\cdots x_{i_k}$, is called a {\em word} over $X$. Its {\em length} is $\abs{\eta}=k$. In particular, $\abs{\eta}_{x_i}$ is the number of times the letter $x_i\in X$ appears in $\eta$. The set of all words including the empty word, $\emptyset$, is denoted by $X^\ast$, and $X^+:=X^\ast\backslash\emptyset$. The set $X^\ast$ forms a monoid under catenation. The set of all words with prefix $\eta$ is written as $\eta X^\ast$. Any mapping $c:X^\ast\rightarrow\re^\ell$ is called a {\em formal power series}. The value of $c$ at $\eta\in X^\ast$ is denoted by $(c,\eta)$ and called the {\em coefficient} of $\eta$ in $c$. A series $c$ is {\em proper} when $(c,\emptyset)=0$ else it is a {\em non-proper} series. The {\em support} of $c$, $\supp(c)$, is the set of all words having nonzero coefficients. The {\em order} of $c$, $\orde(c)$, is the length of the minimal length word in its support. Normally, $c$ is written as a formal sum $c=\sum_{\eta\in X^\ast}(c,\eta)\eta.$ The collection of all formal power series over $X$ is denoted by $\allseriesell$. The $i^{th}$ component of a series $c \in \allseriesell$ is denoted by $c_i$ viz $\left(c_i,\eta\right) = \left(c,\eta\right)_i$. The subset of all proper series in $\allseriesell$ is denoted by $\allproperseries{\ell}$, while the subset of non-proper series is denoted by $\allnproperseries{\ell}$.

\begde A series $c \in \allseries{\ell}$ is called {\em purely improper} if $c_i$ is non-proper $\forall i = 1,\ldots,\ell$. The subset of all purely improper series in $\allseriesell$ is denoted by $\allpiseries{\ell}$.
\endde

Observe that $\allpiseries{\ell} \subsetneq \allnproperseries{\ell}$ if $\ell > 1$, otherwise $\allpiseries{} = \allnproperseries{}$. For the purpose of the document, the product of two vectors in $\re^n$ is given by the Hadamard product. The {\em Cauchy product}, $\mathscr{C}: \allseries{\ell} \times \allseries{\ell} \longrightarrow \allseries{\ell}$ defined as $\left(c,d\right) \mapsto c.d$, where 
\begin{align*}
(c.d,\eta) = \sum_{\substack{\zeta,\nu \in X^{\ast}\\\zeta.\nu = \eta}} \left(c,\zeta\right) 
\left(d,\nu\right)
\end{align*}
Observe that $\allseries{\ell}$ constitutes an associative $\re$-algebra under the Cauchy product. If $d \in \allpiseries{\ell}$, then Cauchy inverse of $d$, denoted by $d^{-1}$ is defined as 
\begin{align*}
\left(d^{-1}\right)_i = \left(d_i,\emptyset\right)^{-1}\left(\sum_{k \in \nat_0} \left(d_i^{\prime}\right)^k\right),
\end{align*}
where $d_i^{\prime} = 1 - \left(d_i/\left(d_i,\emptyset\right)\right)$. Hence, $\allpiseries{\ell}$ forms a group under Cauchy product with $\doubleone = \left[1 1 \cdots 1\right]^{t} \in \re^{\ell}$ as the identity element.
The shuffle product of two words which is a bilinear product uniquely specified by 
\begals
(x_i\eta)\shuffle(x_j\xi)=x_i(\eta\shuffle(x_j\xi))+x_j((x_i\eta)\shuffle \xi),
\endals
where $x_i,x_j\in X$, $\eta,\xi\in X^\ast$ and with $\eta\shuffle\emptyset=\emptyset\shuffle\eta=\eta$ \citep{Fliess_81}. The shuffle product of two series, $\left(c,d\right) \mapsto c\shuffle d$ is defined as
\begin{align*}
(c\shuffle d,\eta) = \sum_{\substack{\zeta,\nu \in X^{\ast}\\\zeta\shuffle\nu = \eta}} \left(c,\zeta\right) \left(d,\nu\right)
\end{align*}
Note that $\allseries{\ell}$ forms an associative and commutative $\re$-algebra under the shuffle product. If $d \in \allpiseries{\ell}$, then shuffle inverse of $d$, denoted by $d^{\shuffle -1}$ is defined as 
\begin{align*}
\left(d^{\shuffle -1}\right)_i = \left(d_i,\emptyset\right)^{-1}\left(\sum_{k \in \nat_0} \left(d_i^{\prime}\right)^{\shuffle k}\right),
\end{align*}
where $d_i^{\prime} = 1 - \left(d_i/\left(d_i,\emptyset\right)\right)$. Hence, $\allpiseries{\ell}$ forms an Abelian group under the shuffle product with $\doubleone = \left[1 1 \cdots 1\right]^{t} \in \re^{\ell}$ as the identity element. The set $\allseriesell$ is an ultrametric space with the ultrametric
\begals
\kappa(c,d) = \sigma^{\orde(c-d)},
\endals
where $c,d \in \allseriesell$ and $\sigma \in \: ]0,1[$. For brevity, $\kappa(c,0)$ is written as $\kappa(c)$, and $\kappa(c,d) = \kappa(c-d)$.
The ultrametric space $(\allseriesell,\kappa)$ is known to be Cauchy complete \citep{Berstel-Reutenauer_88}. The following types of contraction maps will be useful.

\begde Given metric spaces $(E,d)$ and $(E^{\prime},d^{\prime})$, a map $f : E \longrightarrow E^{\prime}$  is said to be a {\em strong
	contraction map} if $\forall s,t \in E$, it satisfies the condition $d^{\prime}(f(s),f(t)) \leq \alpha d(s,t)$ where $\alpha \in [0,1[$.
If $\alpha = 1$, then the map $f$ is said to be a {\em weak contraction map} or a {\em non-expansive map}.
\endde

In the event that the letters of $X$ commute, the set of all formal power series is denoted by $\allcommutingseriesXell$. The formal series with commuting alphabet is indispensable in definition of the formal static maps in Section~\ref{subsec:fromal_static_maps}. For any series $c \in \allcommutingseriesXell$, the natural number $\overline{\omega}(c)$ corresponds to the order of its proper part $c-(c,\emptyset)$.

\subsection{Chen-Fliess Series}

Let $\mathfrak{p}\ge 1$ and $t_0 < t_1$ be given. For a Lebesgue measurable
function $u: [t_0,t_1] \rightarrow\re^m$, define
$\norm{u}_{\mathfrak{p}}=\max\{\norm{u_i}_{\mathfrak{p}}: \ 1 \le i \le m\}$, where $\norm{u_i}_{\mathfrak{p}}$ is the usual
$L_{\mathfrak{p}}$-norm for a measurable real-valued function,
$u_i$, defined on $[t_0,t_1]$.  Let $L^m_{\mathfrak{p}}[t_0,t_1]$
denote the set of all measurable functions defined on $[t_0,t_1]$
having a finite $\norm{\cdot}_{\mathfrak{p}}$ norm and
$B_{\mathfrak{p}}^m(R)[t_0,t_1]:=\{u\in
L_{\mathfrak{p}}^m[t_0,t_1]:\norm{u}_{\mathfrak{p}}\leq R\}$. Given any series $c\in\allseriesell$, the corresponding
{\em Chen-Fliess series} is
\begeq
F_c[u](t) = \sum_{\eta\in X^{\ast}} (c,\eta)\,E_\eta[u](t,t_0), \label{eq:Fliess-operator-defined}
\endeq
where $E_\emptyset[u]=1$ and
\begdi
E_{x_i\bar{\eta}}[u](t,t_0) =
\int_{t_0}^tu_{i}(\tau)E_{\bar{\eta}}[u](\tau,t_0)\,d\tau
\enddi
with $x_i\in X$, $\bar{\eta}\in X^{\ast}$, and $u_0=1$ \citep{Fliess_81}.
If there exists constants $K,M>0$ such that
\begdi
\abs{(c,\eta)}\leq K M^{\abs{\eta}} \abs{\eta}!,\;\; \forall \eta\in X^\ast, \label{eq:Gevrey-growth-condition}
\enddi
then $F_c$ constitutes a well defined mapping from
$B_{\mathfrak p}^m(R)[t_0,$ $t_0+T]$ into $B_{\mathfrak
	q}^{\ell}(S)[t_0, \, t_0+T]$ for sufficiently small $R,T >0$,
where the numbers $\mathfrak{p},\mathfrak{q}\in[1,\infty]$ are
conjugate exponents, i.e., $1/\mathfrak{p}+1/\mathfrak{q}=1$
\citep{Gray-Wang_SCL02}. This map is referred to as a {\em Fliess operator}. Here $\allseriesellLC$ will denote the set of all such {\em locally convergent} generating series. In the absence of any convergence criterion, \rref{eq:Fliess-operator-defined} only defines an operator in a formal sense. 

\subsection{Interconnections of Chen-Fliess series}
Given Chen-Fliess series $F_c$ and $F_d$, where $c,d\in\allseriesell$,
the parallel and product connections satisfy $F_c+F_d=F_{c+d}$ and $F_cF_d=F_{c\shuffle d}$,
respectively\citep{Ree_58,Fliess_81}. The parallel and product connections preserve local convergence and hence the interconnected systems has a Fliess operator representation\citep{Thitsa-Gray_SIAM12,GS_thesis}. When Chen-Fliess series $F_c$ and $F_d$ with $c\in\allseriesXpri{k}$ and $d\in\allseries{\ell}$ are interconnected in a cascade fashion, where $\lvert X^{\prime} \rvert = \ell + 1$, the composite system $F_c\circ F_d$ has a Chen-Fliess series representation $F_{c\circ d}$, where
the {\em composition product} of $c$ and $d$ is given by
\begeq \label{eq:c-circ-d}
c\circ d=\sum_{\eta\in X^{\prime\ast}} (c,\eta)\,\psi_d(\eta)(\mbf{1})
\endeq%
\citep{Ferfera_79,Ferfera_80}. Here $\mbf{1}$ denotes the monomial $1\emptyset$, and $\psi_d$ is the continuous (in the ultrametric sense) algebra homomorphism from $\allseriesXpri{}$ to the set of vector space endomorphisms on $\allseries{}$, $\Endallseries$, uniquely specified by
$\psi_d(x_i^{\prime}\eta)=\psi_d(x_i^{\prime})\circ \psi_d(\eta)$ with
$ 
\psi_d(x_i^{\prime})(e)= x_0(d_i\shuffle e),
$
$i=0,1,\ldots,m$
for any $e\in\allseries{}$, and where $d_i$ is the $i$-th component series of $d$
($d_0:=\mbf{1}$). By definition, $\psi_d(\emptyset)$ is the identity map on $\allseries{}$. The cascade interconnection preserves local convergence and thus the composite has a Fliess operator representation\citep{Thitsa-Gray_SIAM12}. The linearity of the composition product in the left argument is evident from the definition. However, the following theorem states that the shuffle product distributes over the composition product from the left.
\begth\citep{Gray-Li_05}\label{thm:shuff_dist_comp} Let $c,d \in \allseriesXpri{k}$ and $e \in \allseries{\ell}$, such that $\abs{X^{\prime}} = \ell + 1$, then $\left(c\shuffle d\right)\circ e = \left(c \circ e\right) \shuffle \left(d \circ e\right)$. 
\endth
Given a series $e \in \allseries{\ell}$, define a map $\Upsilon_e : \allseriesXpri{k} \longrightarrow \allseries{k}$ defined as $c \mapsto c \circ e$. Theorem~\ref{thm:shuff_dist_comp} infers that $\Upsilon_e$ is an $\re$-algebra homomorphism from the shuffle algebra of $\allseriesXpri{k}$ to the shuffle algebra of $\allseries{\ell}$. The composition product is a strong contraction map with respect to its right argument in the ultrametric topology and is stated in the following theorem.
\begth\citep{Gray-Li_05}\label{thm:comp_strong_cont} Let $c \in \allseriesXpri{k}$ and $d,e \in \allseries{\ell}$, such that $\abs{X^{\prime}} = \ell + 1$, then $\kappa\left(c\circ d, c\circ e\right) \leq \sigma \kappa\left(d,e\right)$.
\endth

\begin{figure}[t]
    	\begin{center}
    		\includegraphics[scale=0.65]{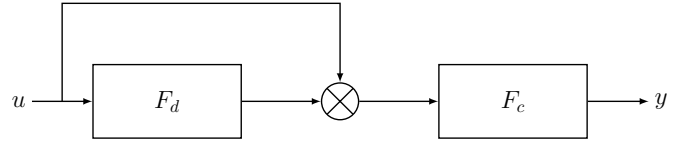}
    	\caption{Cascade connection of Chen-Fliess $F_d$ with $F_c$ along with multiplicative feedforward of input}
    	\label{fig:mult_mix_comp}
    	\end{center}
 \end{figure}

The {\em unital shuffle Chen-Fliess series} arise primarily in the multiplicative output dynamic feedback interconnection of Chen-Fliess series as described in \citet{Gray-KEF_SIAM2017} and Section~\ref{sec:mult_dyn_feedb} of this document. For $\abs{X} = m+1$, the set of all unital shuffle Chen-Fliess series, denoted by $\deltaFliess$, is defined as $\deltaFliess = \{I.F_d : d \in \allseries{m} \}$, where $I$ denotes the identity operator. It is convenient to introduce a symbol $\delta$ as the generating series for the identity map viz. $F_{\delta}[u] = I[u] = u$. Hence, $u. F_{d}[u] = I.F_d[u] = F_{\delta}.F_d[u] = F_{\delta \shuffle d}[u] = F_{\de{d}}[u]$, with $\de{d} = \delta \shuffle d$. The series $\delta \shuffle d$ is the generating series for the Chen-Fliess series depicting the feedforward product of input with the output of $F_d$. The set of all generating series for $\deltaFliess$ shall be denoted by $\delta \shuffle \allseries{m}$. The cascade interconnection of a Chen-Fliess series $F_c$ and $F_{d}$ along with the multiplicative feedforward of the input, as shown in Figure~\ref{fig:mult_mix_comp}, is denoted by $F_{c \mixprocomp \de{d}}$ viz. $F_c[u.F_d[u]] = F_c \circ F_{\de{d}}[u] = F_{c \mixprocomp \de{d}}[u]$, where $c \mixprocomp \de{d}$ denotes the {\em multiplicative mixed composition product} of $c \in \allseries{p}$ and $d \in \allseries{m}$. The multiplicative mixed composition product of $c$ and $d$ , $c \mixprocomp \de{d}$ can be defined as 
\begin{align*}
c \mixprocomp {\de{d}} = \sum_{\eta\in X^\ast} \left(c,\eta\right) \bar{\phi}_d \left(\eta\right)\left(\mbf{1}\right) = \sum_{\eta\in X^\ast} \left(c,\eta\right) \eta \mixprocomp \de{d},
\end{align*}
where $\bar{\phi} : \allseries{} \longrightarrow \Endallseries$ is an $\re$-algebra homomorphism such that $\bar{\phi}_d(x_0)(e) = x_0e$ and $\bar{\phi}_d(x_i)(e) = x_i(d_i \shuffle e)$. Here $\allseries{}$ is taken as an $\re$-algebra under Cauchy product and $\Endallseries$ is an $\re$-algebra under composition. It is straightforward that multiplicative mixed composition product is linear in its left argument. The following results are already known in the single-input single-output (SISO) setting. However, their multi-input multi-output (MIMO) extensions are straightforward and to avoid reiteration of the proofs, only the statements are provided in this document. The foremost of the theorems assert that shuffle product distributes over the multiplicative mixed composition product from the left.  
\begth\citep{Gray-KEF_SIAM2017}\label{thm:shuff_dist_mixprc} Let $c,d \in \allseries{p}$ and $e \in \allseries{m}$, then $\left(c \shuffle d\right) \mixprocomp \de{e} = \left(c \mixprocomp \de{e}\right) \shuffle \\
\left(d \mixprocomp \de{e}\right)$.
\endth

The inference of Theorem~\ref{thm:shuff_dist_mixprc} is that for any $e \in \allseries{m}$, the map $\Gamma_e : \allseries{p} \longrightarrow \allseries{p}$ given by $d \mapsto d \mixcomp \de{e}$ is an $\re$-algebra endomorphism on the shuffle algebra $\allseries{p}$. The next lemma is essential in proving that the multiplicative mixed composition product is a strong contraction map in its right argument in the ultrametric topology.

\begle\citep{Gray-KEF_SIAM2017} Let $\eta \in \allwords$ and $d, e \in \allseries{m}$, then $\kappa\left(\eta \mixprocomp \de{d}, \eta \mixprocomp \de{e}\right) \leq \sigma^{\abs{\eta}} \kappa\left(d,e\right)$. 
\endle  
 
 The following theorem states the strong contraction property of the multiplicative mixed composition product which is an essential result in Section~\ref{sec:mult_dyn_feedb}.
 
 \begth\citep{Gray-KEF_SIAM2017}\label{thm:mixpr_cont} Let $d,e \in \allseries{m}$ and $c \in \allseries{p}$, then $\kappa\left(c \mixprocomp \de{d}, c \mixprocomp \de{e}\right) \leq \sigma^{\ord\left(c^{\prime}\right)}\kappa\left(d,e\right)$, where $c^{\prime}$ is the proper part of $c$  viz. $c^{\prime} = c - \left(c,\emptyset\right)$.
 \endth  
 
 Since $\ord\left(c^{\prime}\right) \geq 1$ and $\sigma \in ]0,1[$, then from Theorem~\ref{thm:mixpr_cont}, the map $\bar{\Gamma}_c : e \mapsto c \mixprocomp \de{e}$ is a strong contraction map in the ultrametric topology. The following lemma is essential in proving the mixed associativity of the composition and multiplicative mixed composition product. The result, along with Theorem~\ref{thm:mix_assoc_comp_mixpr} can be inferred in the SISO setting from Lemma~$3.6$ in \citet{Gray-KEF_SIAM2017}, and its extension to the MIMO case is purely straightforward.

\begle\citep{Gray-KEF_SIAM2017} Let $X^{\prime} = \{x'_0,\ldots,x'_p\}$ and $\eta \in {X^{\prime}}^{\ast}$. Let $d \in \allseries{p}$ and $e \in \allseries{m}$, then $\eta \circ \left(d \mixprocomp \de{de}\right) = \left(\eta \circ d\right)\mixprocomp \de{e}$.
\endle 
 
The following theorem states that the  composition product and multiplicative mixed composition product are associative in combination.

\begth\citep{Gray-KEF_SIAM2017}\label{thm:mix_assoc_comp_mixpr}Let $X^{\prime} = \{x'_0,\ldots,x'_p\}$ and $c \in \allseriesXpri{q}$. Let $d \in \allseries{p}$ and $e \in \allseries{m}$, then $c \circ \left(d \mixprocomp \de{de}\right) = \left(c \circ d\right)\mixprocomp \de{e}$.
\endth

\subsection{Multiplicative Dynamic Output Feedback Group}\label{subsec:dyn_mult_feedb_grp}
The dynamic multiplicative feedback group plays a vital role in computation of the multiplicative dynamic feedback formula, as pictured in SISO setting, in \citet{Gray-KEF_SIAM2017} and in assessing the feedback as a group action in Section~\ref{sec:mult_dyn_feedb}. Consider the cascade interconnection of two unital shuffle Chen-Fliess series $F_{\de{c}}$ and $F_{\de{d}}$, where $c,d \in \allseries{m}$. The composite system is given by the Chen-Fliess series $F_{\de{c}\circ\de{d}}$, where $\de{c}\circ \de{d}$ denotes the {\em multiplicative composition product} of $\de{c}$ and $\de{d}$ and is defined as 
\begin{align}\label{eqn:grp_prod}
\de{c}\circ \de{d} = \delta \shuffle \left(d \shuffle c \mixprocomp \de{d}\right) = \de{\left(d \shuffle c \mixprocomp \de{d}\right)}.
\end{align}
There is an abuse of notation $\circ$ between \rref{eq:c-circ-d} and \rref{eqn:grp_prod}, however the meaning of $\circ$ should always be clear from the context. The following theorem states the multiplicative composition product is associative. The result, along with Theorem~\ref{thm:mixprocomp_right_act_monoid} were stated and proven in Lemma~$3.6$ of \citet{Gray-KEF_SIAM2017} in the SISO setting but the authors' proofs were independent of the SISO. Hence, the statements along with the proofs naturally extend to the MIMO setting. 
\begth\citep{Gray-KEF_SIAM2017}\label{thm:mult_comp_assoc} Let $c,d,e \in \allseries{m}$, then, $\left(\de{c}\circ\de{d}\right)\circ \de{e} = \de{c}\circ\left(\de{d}\circ \de{e}\right)$. 
\endth
Observe that \rref{eqn:grp_prod} and Theorem~\ref{thm:mult_comp_assoc} infer that $\delta \shuffle \allseries{m}$ forms a noncommutative monoid under multiplicative composition product, with the identity element $\de{\doubleone}$. The following theorem states that the multiplicative mixed composition product is right action on $\allseries{q}$ by the monoid $\left(\delta \shuffle \allseries{m},\circ\right)$.
\begth\citep{Gray-KEF_SIAM2017}\label{thm:mixprocomp_right_act_monoid} Let $c \in \allseries{q}$ and $d,e \in \allseries{m}$, then $\left(c \mixprocomp \de{d}\right)\mixprocomp \de{e} = c \mixprocomp \left(\de{d} \circ \de{e}\right)$.
\endth
The prominent question is to find the invertible elements of the monoid $\left(\delta \shuffle \allseries{m},\circ\right)$. Let $d,e \in \allpiseries{m}$ and suppose
\begin{align*}
\de{d} \circ \de{e} &= \de{\doubleone}
\end{align*}
Applying~\rref{eqn:grp_prod},
\begin{align*}
e \shuffle \left(d \mixprocomp \de{e}\right) = \doubleone 
\end{align*}
Observe that $d \in \allpiseries{m}$ implies $\d \mixprocomp \de{e} \in \allpiseries{m}$ and using Theorem~\ref{thm:shuff_dist_mixprc},
\begin{align*}
e = \left(d \mixprocomp \de{e}\right)^{\shuffle -1} = d^{\shuffle -1} \mixprocomp \de{e}.
\end{align*}
Hence, for $\de{e}$ to be right inverse of $\de{d}$, the purely improper series $e$ has to satisfy the fixed point equation 
\begin{align}\label{eqn:fixed_point_grp_eqn_rinv}
e = d^{\shuffle -1} \mixprocomp \de{e}
\end{align}
Observe from Theorem~\ref{thm:mixpr_cont} that the map $e \mapsto d^{\shuffle -1} \mixprocomp \de{e}$ is a strong contraction in the ultrametric space inferring that \rref{eqn:fixed_point_grp_eqn_rinv} has a unique fixed point. Suppose $\de{e}$ is the left inverse of $\de{d}$ viz $\de{e} \circ \de{d}$, then a similar procedure shows that $e$ has to satisfy the equation
\begin{align}\label{eqn:grp_eqn_linv}
d = e^{\shuffle -1} \mixprocomp \de{d}
\end{align}
Note that if $e$ is a solution of \rref{eqn:fixed_point_grp_eqn_rinv}, then $e$ satisfies \rref{eqn:grp_eqn_linv} and converse. Hence, $e$ is given the notation $d^{\circ -1}$ and for $d \in \allpiseries{m}$, the inverse of $\de{d}$ exists and is unique, denoted by $\de{d^{\circ -1}}$ viz. $\left(\de{d}\right)^{\circ -1} = \delta \shuffle d^{\circ -1} = \de{d^{\circ -1}}$. Thus, $\delta \shuffle \allpiseries{m}$ forms a group under multiplicative composition product, termed as the {\em multiplicative dynamic output feedback group} and is formally stated in the following theorem.
\begth\label{thm:mult_dyn_feedb_grp} $\left(\delta \shuffle \allpiseries{m},\circ\right)$ forms a group with the identity element $\de{1}$ 
\endth
It is worth noting that \citet{Gray-KEF_SIAM2017} proved Theorem~\ref{thm:mult_dyn_feedb_grp} for one-dimensional case viz. $m=1$. In light of Theorem~\ref{thm:mult_dyn_feedb_grp}, Theorem~\ref{thm:shuff_dist_mixprc} and \rref{eqn:grp_prod} one obtains the following relations for $c \in \allpiseries{m}$:
\begin{align}\label{eqn:inverse_relations}
c^{\circ -1} &= c^{\shuffle -1} \mixprocomp \de{c^{\circ -1}}\\
\left(c^{\circ -1}\right)^{\shuffle -1} &= c \mixprocomp \de{c^{\circ -1}}
\end{align} 

\subsection{Cauchy Algebra of Formal Static Maps}\label{subsec:fromal_static_maps}
This subsection provides a brief discussion on formal static maps, which are used to describe the memoryless maps in the feedback path of static feedback interconnection, as described in Section~\ref{sec:mult_stat_feedb}. Let $\tilde{X} = \{\tilde{x}_1,\ldots,\tilde{x}_m\}$ and $d \in \allcommutingseriesX{k}$. A formal static function $f_d : \re^m \longrightarrow \re^k$  around the point $z = 0$ is defined as 
    \begin{align*}
    f_d\left(z\right) &= \sum_{\eta \in \tilde{X}^{\ast}} \left(d,\eta\right) z^{\eta},
    \end{align*}
    where $z \in \re^m$, and $z^{\tilde{x}_i\eta}= z_iz^{\eta} \; \forall \tilde{x}_i \in \tilde{X}, \eta\in \allcommutingwords$. The base case is taken to be $z^{\emptyset} = 1$. Denote the collection of all formal static maps from $\re^m$ to $\re^k$ as $\Homstat{m}{k}$. The series $d \in \allcommutingseriesX{k}$ is called the generating series of the static map $f_d$. A series $d \in \allcommutingseriesX{}$ is said to be \emph{locally convergent} if there exist constants $K_d,M_d >0$ such that $\abs{\left(d,\eta\right)} \leq K_dM_d^{\abs{\eta}}, \; \forall \eta \in \tilde{X}^{\ast}$. A series $d \in \allcommutingseriesX{k}$ is said to be locally convergent if and only if each component $d_i$ is locally convergent for $i=1,\ldots, m$. The subset of all locally convergent series in $\allcommutingseriesX{k}$ is denoted as $\allcommutingseriesLCX{k}$. The following theorem describes the significance of the definition of local convergence in the context of formal static maps. The condition $d \in \allcommutingseriesLCX{}$ is both necessary and sufficient for the corresponding static function $f_d$ to be a locally (around $z=0$) real analytic map\citep{GS_thesis}. The following states the formal static maps are closed under pointwise multiplication in $\re^k$.
\begth\label{thm:cauchy_inv_static} Let the formal static maps $f_d,f_e : \re^m \longrightarrow \re^k$, with $d,e \in \allcommutingseriesX{k}$. The product of the maps $f_d.f_e : \re^m \longrightarrow \re^k$ is a formal static map $f_{d.e}$, where $d.e$ is the Cauchy product of $d$ and $e$.
\endth
\begpr
\begin{align*}
f_d.f_e(z) &= f_d(z).f_e(z) = \left(\sum_{\zeta \in \tilde{X}^{\ast}} \left(d,\zeta\right) z^{\zeta}\right)\left(\sum_{\nu \in \tilde{X}^{\ast}} \left(e,\nu\right) z^{\nu}\right)\\
&= \left(\sum_{\eta \in \tilde{X}^{\ast}} \sum_{\substack{\zeta,\nu \in \tilde{X}^{\ast}\\ \zeta\nu = \eta}}\left(d,\zeta\right)\left(e,\nu\right) z^{\zeta}z^{\nu}\right)\\
&= \sum_{\eta \in \tilde{X}^{\ast}} \left(d.e,\eta\right) z^{\eta} = f_{d.e}
\end{align*}
\endpr

Theorem~\ref{thm:cauchy_inv_static} asserts that $\Homstat{m}{k}$ forms an $\re$-algebra there is an $\re$-algebra isomorphism from the $\Homstat{m}{k}$ to the Cauchy algebra of $\allcommutingseries{k}$. Let $f_d$ be  formal static map, with $d \in \allcommutingpiseries{k}$. Then from Theorem~\ref{thm:cauchy_inv_static}, the generating series of the multiplicative inverse of the formal static map $f_d$, denoted by $f_d^{-1}$, is given by the Cauchy iverse of $d$ viz. $f_d^{-1} = f_{d^{-1}}$. Hence, the unit group of $\Homstat{m}{k}$ is isomorphic to the group $\allcommutingpiseries{k}$ under Cauchy product. 

\subsection{Wiener-Fliess Composition of Formal Power Series}\label{subsec:Wiener-Fliess}
\begin{figure}[t]
    	\begin{center}
    		\includegraphics[scale = 0.35]{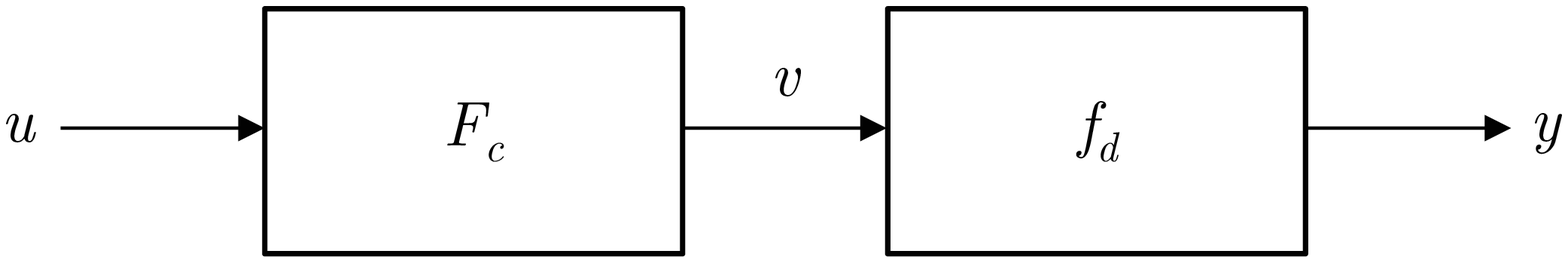}
    	\caption{Wiener-Fliess connection}
    	\label{fig:wiener-fliess}
    	\end{center}
 \end{figure}

This subsection describes the cascade connection shown in Figure~\ref{fig:wiener-fliess} of a Chen-Fliess series $F_c$ generated by a series $c \in \allproperseriesell$ and a formal static map $f_d \in \Homstat{\ell}{k}$. Such configurations are called \emph{Wiener-Fliess} connections. The connection is known to generate well-defined Chen-Fliess series for the composite system, and its generating series is computed through the \emph{Wiener-Fliess composition product}. The definition of Wiener-Fliess composition product first appeared in \citet{Thitsa-Gray_SIAM12}, however, the definition was expanded even for $c \in \allnproperseries{\ell}$ in \citet{GS_thesis}. However, the current document works with the restricted definition.   
    
\begth\citep{Gray_Thitsa_IJC12,Venkat-Gray_2021} \label{thm:Wiener-Fliess-product}
Let $ X = \{x_0,x_1,\ldots,$ $x_m\}$ and $\tilde{X} = \{\tilde{x}_1,\tilde{x}_2,\ldots,\tilde{x}_\ell\}$. Given a formal Fliess operator $F_c$ with $c \in \allproperseriesell$ and formal function $f_d \in \Homstat{\ell}{k}$, then the composition $f_d\circ F_c$ has a generating series in $\allseries{k}$ given by the Wiener-Fliess composition product
    \begin{align}
    d \mixcomp c = \sum_{\tilde{\eta} \in \tilde{X}^{\ast}} (d,\tilde{\eta}) c^{\shuffle \tilde{\eta}}\label{eqn:mixcomp},
    \end{align}
    where $c^{\shuffle\,\tilde{x}_{i}\tilde{\eta}}:=c_{i}\shuffle c^{\shuffle \tilde{\eta}} \; \forall \tilde{x}_i \in \tilde{X} , \; \forall \tilde{\eta} \in \allcommutingwords$, and $c^{\shuffle \phi} = 1$.
    \endth  
    
 Observe that if $d \in \allcommutingpiseries{k}$ and $c \in \allproperseries{\ell}$, then $d\mixcomp c \in \allpiseries{k}$. For a fixed $d \in \allcommutingseriesX{k}$ define the map $d_{\mixcomp} : \allproperseriesell \longrightarrow \allseries{k}:c \mapsto d\mixcomp c$. The Wiener-Fliess connection preserves local convergence and hence, the composite system has a Fliess operator representation\citep{GS_thesis}. The following theorems describe the contractive properties of $d_{\mixcomp}$ and $\tilde{d}_{\mixcomp}$ in the ultrametric topology.    
 \begth\citep{Venkat-Gray_2021} \label{th:Wiener-Fliess-contraction}The map $d_{\mixcomp}$ is a weak contraction map when $\overline{\omega}\:(d) = 1$ and a strong contraction map when $\overline{\omega}\:(d) > 1$.
 \endth     
Theorem~\ref{th:Wiener-Fliess-contraction} infers that the Wiener-Fliess composition product is at the very least is a weak contraction map with respect to the noncommutative formal series argument. The following theorem is a crucial result used in computing the feedback formula for the multiplicative static feedback connection.

\begth\label{thm:Cauchy_prod_WF} Let $d,d^{\prime} \in \allcommutingseriesX{k}$ and $c \in \allproperseries{\ell}$, then $\left(d.d^{\prime}\right) \mixcomp c = \left(d\mixcomp c\right) \shuffle \left(d^{\prime}\mixcomp c\right)$ provided $d\mixcomp c$ and $d^{\prime} \mixcomp c$ are well-defined.
\endth
\begpr
\begin{align*}
\left(d.d^{\prime}\right) \mixcomp c &= \sum_{\eta \in \tilde{X}^{\ast}}\left(d.d^{\prime},\eta\right) c^{\shuffle \eta} \\
&= \sum_{\eta \in \tilde{X}^{\ast}} \sum_{\substack{\zeta,\nu \in \tilde{X}^{\ast}\\\zeta\nu = \eta}} \left(d,\zeta\right)\left(d^{\prime},\nu\right) z^{\shuffle \zeta}\shuffle z^{\shuffle \nu} \\
&= \left(\sum_{\zeta \in \tilde{X}^{\ast}}\left(d,\zeta\right) z^{\shuffle \zeta}\right) \shuffle \left(\sum_{\nu \in \tilde{X}^{\ast}}\left(d^{\prime},\nu\right) z^{\shuffle \nu}\right) \\
&= \left(d \mixcomp c\right) \shuffle \left(d^{\prime} \mixcomp c\right). 
\end{align*}
\endpr
  
Note that $(\allcommutingseriesX{k},\cdot)$ is a commutative monoid as the letters of $\tilde{X}$ commute. For a series $c \in \allproperseriesell$, define a map $\Omega_c : \allcommutingseriesX{k} \longrightarrow \allseries{k}$. Theorem~\ref{thm:Cauchy_prod_WF} asserts that $\Omega_c$ is an $\re$-algebra homomorphism from the Cauchy algebra of $\allcommutingseriesX{k}$ to the shuffle algebra of $\allseries{k}$. If $c$ is proper, then via theorem~\ref{thm:Cauchy_prod_WF}, it is evident that $d^{-1}\mixcomp c = \left(d \mixcomp c\right)^{\shuffle -1}$, provided $d \in \allcommutingpiseries{k}$. Hence, there is a group homomorphism from the $\allcommutingpiseries{k}$ group under Cauchy product to the $\allpiseries{k}$ group under shuffle product via the map $\Omega_c$. The following theorem is pivotal in Section~\ref{sec:mult_stat_feedb} and states Wiener-Fliess compositon product and multiplicative mixed composition product are mixed associative.
  
\begth\label{thm:WF_mixprocomp_mixassoc} Let $d \in \allcommutingseriesX{p}, c \in \allseries{q}$ and $e \in \allseries{m}$, such that $\abs{\tilde{X}} = q$, then $d \mixcomp \left(c \mixprocomp \de{e}\right) = \left(d \mixcomp c\right) \mixprocomp \de{e}$.
\endth

\begpr
\begin{align*}
d \mixcomp \left(c \mixprocomp \de{e}\right) = \sum_{\eta \in \tilde{X}^{\ast}} \left(d,\eta\right) \left(c \mixprocomp \de{e}\right)^{\shuffle \eta}.
\end{align*}
Using Theorem~\ref{thm:shuff_dist_mixprc},
\begin{align*}
d \mixcomp \left(c \mixprocomp \de{e}\right) &= \sum_{\eta \in \tilde{X}^{\ast}} \left(d,\eta\right) \left(c^{\shuffle \eta} \mixprocomp \de{e}\right) \\
&= \left[\sum_{\eta \in \tilde{X}^{\ast}} \left(d,\eta\right) c^{\shuffle \eta}\right] \mixprocomp \de{e}\\ 
&=  \left(d \mixcomp c\right) \mixprocomp \de{e}.
\end{align*}
\endpr

\section{Chen-Fliess Series Under Multiplicative Dynamic Output Feedback}\label{sec:mult_dyn_feedb}
\begin{figure}[tb]
	\vspace*{-0.4 in}
    	\begin{center}
    		\includegraphics[scale = 0.7]{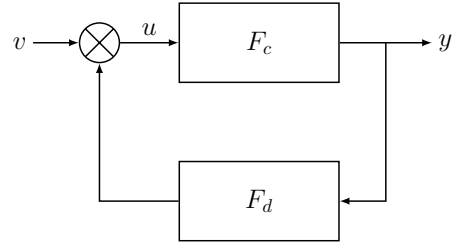}
    	\caption{Chen-Fliess series $F_c$ in multiplicative output feedback with Chen-Flies series $F_d$}
    	\label{fig:mult_dyn_out_feedb}
    	\end{center}
 \end{figure}

Let $F_c$ be a Chen-Fliess series with a generating series $c \in \allseries{q}$. Assume it is interconnected with a Chen-Fliess series $F_d$ with a purely improper generating series $d \in \allpiseriesXpri{m}$, as shown in Figure~\ref{fig:mult_dyn_out_feedb}. Note that, $\abs{X} = m+1$ and $\abs{X^{\prime}} = q+1$. The primary goal of this section is to show that the closed-loop system has a Chen-Fliess series representation, say $y = F_e[v]$, where $e \in \allseries{q}$. If this is the case, then necessarily 
\begin{align*}
y &= F_e[v] = F_c[u] = F_c[vF_d[y]] \\
&= F_c[vF_d[F_e[v]]] = F_c[vF_{d\circ e}[v]]\\
&= F_{c \mixprocomp \de{\left(d \circ e\right)}}[v]
\end{align*}   
for any admissible input $u$. Therefore, the series $e$ has to satisy the fixed point equation
\begin{align}\label{eqn:dyn_feedb_fixed_point_eqn}
e = c \mixprocomp \de{\left(d\circ e\right)}.
\end{align}
Observe that, in light of Theorem~\ref{thm:comp_strong_cont} and Theorem~\ref{thm:mixpr_cont} the map $e \mapsto c \mixprocomp \de{\left(d \circ e \right)}$ is a strong contraction map in the ultrametric space and thus \rref{eqn:dyn_feedb_fixed_point_eqn} has a unique fixed point. The following thoerem establishes the first main result of this section, which follows immediately.  
\begth\label{thm:fixed_point_dyn_feedb} The series $c \mixprocomp \de{\left(d^{\shuffle -1}\circ c\right)^{\circ -1}} \in \allseries{q}$ is the unique fixed point of the map  $e \mapsto c \mixprocomp \de{\left(d \circ e \right)}$.  
\endth
\begpr If $e := c \mixprocomp \de{\left(d^{\shuffle -1}\circ c\right)^{\circ -1}}$, then
\begin{align*}
 c \mixprocomp \de{\left(d \circ e \right)} &= c \mixprocomp \de{\left[d \circ \left(c \mixprocomp \de{\left(d^{\shuffle -1}\circ c\right)^{\circ -1}}\right)\right]} \\
\end{align*}
Using Theorem~\ref{thm:mix_assoc_comp_mixpr} and then Theorem~\ref{thm:shuff_dist_mixprc}, 
\begin{align*}
c \mixprocomp \de{\left(d \circ e \right)} &= c \mixprocomp \de{\left[\left(d\circ c\right) \mixprocomp \de{\left(d^{\shuffle -1}\circ c\right)^{\circ -1}} \right]}\\
&= c \mixprocomp \de{\left[\left(d \circ c \right)^{\shuffle -1} \mixprocomp \de{\left(d^{\shuffle-1}\circ c\right)^{\circ -1}}\right]^{\shuffle -1}}.
\end{align*}
Using Theorem~\ref{thm:shuff_dist_comp},
\begin{align*}
c \mixprocomp \de{\left(d \circ e \right)} &= c \mixprocomp \de{\left[\left(d^{\shuffle -1}\circ c\right) \mixprocomp \de{\left(d^{\shuffle-1}\circ c\right)^{\circ -1}}\right]}^{\shuffle -1}.
\end{align*}
Using the relations~\rref{eqn:inverse_relations},
\begin{align*}
c \mixprocomp \de{\left(d \circ e \right)} &= c \mixprocomp \de{\left[\left(\left(d^{\shuffle -1}\circ c\right)^{\circ -1}\right)^{\shuffle -1}\right]^{\shuffle -1}} \\
&= c \mixprocomp \de{\left(d^{\shuffle -1} \circ c\right)^{\circ -1}} = e.
\end{align*}
\endpr
\begth\label{thm:mult_dyn_feedb_formula} Given a series $c \in \allseries{q}$ and a purely improper series $d \in \allpiseriesXpri{m}$ (such that $\abs{X} = m+1$ and $\abs{X^{\prime}} = q+1$), then the generating series for the closed-loop system in Figure is given by the {\em multiplicative dynamic feedback product} $c\check{@}d := c \mixprocomp \de{\left(d^{\shuffle -1}\circ c\right)^{\circ -1}}$.
\endth

The notion that feedback can described mathematically as a transformation group acting on the plant is well established in control theory\citep{Brockett_78}. The following theorem describes the situation in the present context.

\begth\label{thm:mult_dyn_feedb_grp_action} The multiplicative dynamic feedback product is a right group action by the multiplicative group $\left(\allpiseriesXpri{m},\shuffle,\doubleone\right)$ on the set $\allseries{q}$, where $\abs{X} = m+1$ and $\abs{X^{\prime}} = q+1$.  
\endth
\begpr Let $c \in \allseries{q}$. Observe that from Theorem~\ref{thm:mult_dyn_feedb_formula},
\begin{align*}
c \check{@}\doubleone &= c\mixprocomp \de{\left(\doubleone^{\shuffle -1} \circ c\right)^{\circ -1} } \\
&= c \mixprocomp \de{\doubleone} = c.
\end{align*}
Let $d_1,d_2 \in \allpiseriesXpri{m}$. It needs to be proven that $\left(c\check{@}d_1\right)\check{@}d_2 = c\check{@}\left(d_1 \shuffle d_2\right)$. 
From Theorem~\ref{thm:mult_dyn_feedb_formula}, observe that
\begin{align*}
\left(c\check{@}d_1\right)\check{@}d_2 &= \left(c\check{@}d_1\right) \mixprocomp \de{\left(d_2^{\shuffle -1}\circ \left(c \check{@}d_1\right)\right)^{\circ -1}} \\
&= \left(c \mixprocomp \de{\left(d_1^{\shuffle -1} \circ c\right)^{\circ -1}}\right) \mixprocomp \\
& \quad \quad \de{\left(d_2^{\shuffle -1}\circ \left(c \mixprocomp \de{\left(d_1^{\shuffle -1}\circ c\right)^{\circ -1}}\right)\right)^{\circ -1}}.
\end{align*}
Applying Theorem~\ref{thm:mix_assoc_comp_mixpr},
\begin{align*}
\left(c\check{@}d_1\right)\check{@}d_2 &= \left(c \mixprocomp \de{\left(d_1^{\shuffle -1} \circ c\right)^{\circ -1}}\right) \mixprocomp \\
& \quad \quad \de{\left(\left(d_2^{\shuffle -1}\circ c\right)\mixprocomp \de{\left(d_1^{\shuffle -1}\circ c\right)^{\circ -1}}\right)^{\circ -1}}.
\end{align*}
Applying Theorem~\ref{thm:mixprocomp_right_act_monoid} and fact that the group inverse obeys $\de{g}^{\circ -1} \circ \de{h}^{\circ -1} = \left(\de{h}\circ \de{g}\right)^{\circ -1} \, \forall \; \de{g},\de{h} \in \delta \shuffle \allpiseries{m}$
\begin{align*}
\left(c\check{@}d_1\right)\check{@}d_2 &= c \mixprocomp \Big[ \de{\left(d_1^{\shuffle -1}\circ c\right)^{\circ -1}} \circ \\
&\quad \quad \de{\left(\left(d_2^{\shuffle -1}\circ c\right)\mixprocomp \de{\left(d_1^{\shuffle -1}\circ c\right)^{\circ -1}}\right)^{\circ -1}} \Big] \\
&=c \mixprocomp \Big[ \de{\left(\left(d_2^{\shuffle -1}\circ c\right) \mixprocomp \de{\left(d_1^{\shuffle -1}\circ c\right)^{\circ -1}}\right)} \\
&\quad \quad \quad \circ \de{\left(d_1^{\shuffle -1}\circ c\right)}\Big]^{\circ -1}.
\end{align*}
Applying~\rref{eqn:grp_prod},
\begin{align*}
\left(c\check{@}d_1\right)\check{@}d_2 &= c \mixprocomp \de{}\Bigg[\left(d_1^{\shuffle -1} \circ c\right) \shuffle
\Bigg(\Big(\left(d_{2}^{\shuffle -1}\circ c\right) \mixprocomp \\
&\quad \quad \de{\left(d_{1}^{\shuffle -1}\circ c\right)^{\circ -1}}\Big) \mixprocomp \de{\left(d_1^{\shuffle -1}\circ c\right)}\Bigg)\Bigg]^{\circ -1}.
\end{align*}
Using Theorem~\ref{thm:mixprocomp_right_act_monoid},
\begin{align*}
\left(c\check{@}d_1\right)\check{@}d_2 &= c \mixprocomp \de{}\Bigg[\left(d_1^{\shuffle -1} \circ c\right) \shuffle
\Bigg(\left(d_{2}^{\shuffle -1}\circ c\right) \mixprocomp \\
& \quad \quad \left(\de{\left(d_1^{\shuffle -1}\circ c\right)^{\circ -1}} \circ \de{\left(d_1^{\shuffle -1}\circ c\right)}\right)\Bigg)\Bigg]^{\circ -1} \\
&= c \mixprocomp \de{\left(\left(d_{1}^{\shuffle -1}\circ c\right) \shuffle \left(\left(d_2^{\shuffle -1} \circ c\right) \mixprocomp \de{\doubleone}\right)\right)^{\circ -1}} \\
&= c \mixprocomp \de{\left(\left(d_{1}^{\shuffle -1}\circ c\right) \shuffle \left(d_{2}^{\shuffle -1}\circ c\right)\right)^{\circ -1}}.
\end{align*}
In light of Theorem~\ref{thm:shuff_dist_comp},
\begin{align*}
\left(c\check{@}d_1\right)\check{@}d_2 &= c \mixprocomp \de{\left(\left(d_1^{\shuffle -1}\shuffle d_2^{\shuffle -1}\right)\circ c\right)^{\circ -1}} \\
&= c \mixprocomp \de{\left(\left(d_1\shuffle d_2\right)^{\shuffle -1}\circ c\right)^{\circ -1}} \\
&= c \check{@}\left(d_1 \shuffle d_2\right).
\end{align*}
\endpr

It is worth noting that for the {\em additive dynamic feedback product} the transformation group is the additive group $(\allseriesXpri{m},+,0)$ while here $(\allpiseriesXpri{m},\shuffle,\doubleone)$ plays the role. 

\section{Chen-Fliess Series Under Multiplicative Static Output Feedback}\label{sec:mult_stat_feedb}
\begin{figure}[tb]
	\vspace*{-0.4in}
    	\begin{center}
    		\includegraphics[scale = 0.7]{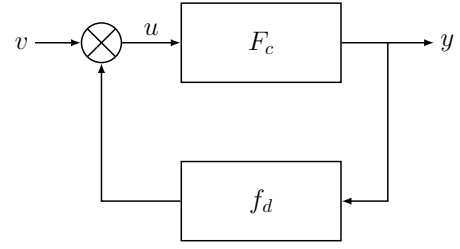}
    	\caption{Chen-Fliess series $F_c$ in multiplicative output feedback with $f_d$}
    	\label{fig:mult_stat_out_feedb}
    	\end{center}
 \end{figure}

Let $F_c$ be a Chen-Fliess series with a proper generating series $c \in \allproperseries{q}$. Assume it is interconnected with a formal static map $f_d$ with a purely improper generating series $d \in \allcommutingpiseries{m}$, as shown in Figure~\ref{fig:mult_stat_out_feedb}. Note that, $\abs{X} = m+1$ and $\abs{\tilde{X}} = q$. The primary goal of this section is to show that the closed-loop system has a Chen-Fliess series representation, say $y = F_e[v]$, where $e \in \allseries{q}$. If this is the case, then necessarily 
\begin{align*}
y &= F_e[v] = F_c[u] = F_c[vf_d[y]] \\
&= F_c[vf_d[F_e[v]]] = F_c[vF_{d\mixcomp e}[v]]\\
&= F_{c \mixprocomp \de{\left(d \mixcomp e\right)}}[v]
\end{align*}   
for any admissible input $u$. Therefore, the series $e$ has to satisy the fixed point equation
\begin{align}\label{eqn:stat_feedb_fixed_point_eqn}
e = c \mixprocomp \de{\left(d\mixcomp e\right)}.
\end{align}
In addition, $e$ must be a proper series for the Wiener-Fliess composition $d\mixcomp e$ to be well defined for arbitrary $d \in \allcommutingpiseries{m}$. It follows directly from the definition of the multiplicative
mixed composition product that if $c \in \allproperseries{\ell}$ then $c\mixprocomp \de{w}$ is also a proper
series for all $w \in \allseries{m}$.

Observe that, in light of Theorem~\ref{th:Wiener-Fliess-contraction} and Theorem~\ref{thm:mixpr_cont} the map $e \mapsto c \mixprocomp \de{\left(d \mixcomp e \right)}$ is a strong contraction map in the ultrametric space and thus \rref{eqn:stat_feedb_fixed_point_eqn} has a unique fixed point. The following fixed point theorem establishes the first main result of this section, which follows immediately.  

\begth\label{thm:fixed_point_stat_feedb} The series $c \mixprocomp \de{\left(d^{-1}\mixcomp c\right)^{\circ -1}} \in \allproperseries{q}$ is the unique fixed point of the map  $e \mapsto c \mixprocomp \de{\left(d \mixcomp e \right)}$.  
\endth
\begpr If $e := c \mixprocomp \de{\left(d^{-1}\mixcomp c\right)^{\circ -1}}$, then
\begin{align*}
 c \mixprocomp \de{\left(d \mixcomp e \right)} &= c \mixprocomp \de{\left[d \mixcomp \left(c \mixprocomp \de{\left(d^{-1}\mixcomp c\right)^{\circ -1}}\right)\right]} \\
\end{align*}
Using Theorem~\ref{thm:WF_mixprocomp_mixassoc} and then Theorem~\ref{thm:shuff_dist_mixprc}, 
\begin{align*}
c \mixprocomp \de{\left(d \mixcomp e \right)} &= c \mixprocomp \de{\left[\left(d\mixcomp c\right) \mixprocomp \de{\left(d^{-1}\mixcomp c\right)^{\circ -1}} \right]}\\
&= c \mixprocomp \de{\left[\left(d \mixcomp c \right)^{\shuffle-1} \mixprocomp \de{\left(d^{-1}\mixcomp c\right)^{\circ -1}}\right]^{\shuffle -1}}.
\end{align*}
Using Theorem~\ref{thm:Cauchy_prod_WF},
\begin{align*}
c \mixprocomp \de{\left(d \mixcomp e \right)} &= c \mixprocomp \de{\left[\left(d^{-1}\mixcomp c\right) \mixprocomp \de{\left(d^{-1}\mixcomp c\right)^{\circ -1}}\right]^{\shuffle -1}}.
\end{align*}
Using the relations~\rref{eqn:inverse_relations},
\begin{align*}
c \mixprocomp \de{\left(d \mixcomp e \right)} &= c \mixprocomp \de{\left[\left(\left(d^{-1}\mixcomp c\right)^{\circ -1}\right)^{\shuffle -1}\right]^{\shuffle -1}} \\
&= c \mixprocomp \de{\left(d^{-1} \mixcomp c\right)^{\circ -1}} = e.
\end{align*}
\endpr
\begth\label{thm:mult_stat_feedb_formula} Given a series $c \in \allproperseries{q}$ and a purely improper series $d \in \allcommutingpiseries{m}$ (such that $\abs{X} = m+1$ and $\abs{\tilde{X}} = q$), then the generating series for the closed-loop system in Figure is given by the {\em multiplicative static feedback product} $c\bar{@}d := c \mixprocomp \de{\left(d^{-1}\circ c\right)^{\circ -1}}$.
\endth

The following theorem describes the transformation group on the plant which characterizes the multiplicative static feedback product. 
\begth\label{thm:mult_stat_feedb_grp_action} The multiplicative static feedback product is a right group action by the Abelian multiplicative group $\left(\allcommutingpiseries{m},\cdot,\doubleone\right)$ on the set $\allproperseries{q}$, where $\abs{X} = m+1$ and $\abs{\tilde{X}} = q$.  
\endth
\begpr Let $c \in \allseries{q}$. Observe that from Theorem~\ref{thm:mult_stat_feedb_formula},
\begin{align*}
c \bar{@}\doubleone &= c\mixprocomp \de{\left(\doubleone^{-1} \mixcomp c\right)^{\circ -1} } \\
&= c \mixprocomp \de{\doubleone} = c.
\end{align*}
Let $d_1,d_2 \in \allcommutingpiseries{m}$. It needs to be proven that $\left(c\bar{@}d_1\right)\bar{@}d_2 = c\bar{@}\left(d_1.d_2\right)$.
From Theorem~\ref{thm:mult_stat_feedb_formula}, observe that
\begin{align*}
\left(c\bar{@}d_1\right)\bar{@}d_2 &= \left(c\bar{@}d_1\right) \mixprocomp \de{\left(d_2^{-1}\mixcomp \left(c \bar{@}d_1\right)\right)^{\circ -1}} \\
&= \left(c \mixprocomp \de{\left(d_1^{-1} \mixcomp c\right)^{\circ -1}}\right) \mixprocomp \\
& \quad \quad \de{\left(d_2^{-1}\mixcomp \left(c \mixprocomp \de{\left(d_1^{-1}\mixcomp c\right)^{\circ -1}}\right)\right)^{\circ -1}}.
\end{align*}
Applying Theorem~\ref{thm:WF_mixprocomp_mixassoc},
\begin{align*}
\left(c\bar{@}d_1\right)\bar{@}d_2 &= \left(c \mixprocomp \de{\left(d_1^{-1} \mixcomp c\right)^{\circ -1}}\right) \mixprocomp \\
& \quad \quad \de{\left(\left(d_2^{-1}\mixcomp c\right)\mixprocomp \de{\left(d_1^{-1}\mixcomp c\right)^{\circ -1}}\right)^{\circ -1}}.
\end{align*}
Applying Theorem~\ref{thm:mixprocomp_right_act_monoid} and fact that the group inverse obeys $\de{g}^{\circ -1} \circ \de{h}^{\circ -1} = \left(\de{h}\circ \de{g}\right)^{\circ -1} \, \forall \; \de{g},\de{h} \in \delta \shuffle \allpiseries{m}$
\begin{align*}
\left(c\bar{@}d_1\right)\bar{@}d_2 &= c \mixprocomp \Big[ \de{\left(d_1^{-1}\mixcomp c\right)^{\circ -1}} \circ \\
&\quad \quad \de{\left(\left(d_2^{-1}\mixcomp c\right)\mixprocomp \de{\left(d_1^{-1}\mixcomp c\right)^{\circ -1}}\right)^{\circ -1}} \Big] \\
&=c \mixprocomp \Big[ \de{\left(\left(d_2^{-1}\mixcomp c\right) \mixprocomp \de{\left(d_1^{-1}\mixcomp c\right)^{\circ -1}}\right)} \\
&\quad \quad \quad \circ \de{\left(d_1^{-1}\mixcomp c\right)}\Big]^{\circ -1}.
\end{align*}
Applying~\rref{eqn:grp_prod},
\begin{align*}
\left(c\bar{@}d_1\right)\bar{@}d_2 &= c \mixprocomp \de{}\Bigg[\left(d_1^{-1} \mixcomp c\right) \shuffle
\Bigg(\Big(\left(d_{2}^{-1}\mixcomp c\right) \mixprocomp \\
&\quad \quad \de{\left(d_{1}^{-1}\mixcomp c\right)^{\circ -1}}\Big) \mixprocomp \de{\left(d_1^{-1}\mixcomp c\right)}\Bigg)\Bigg]^{\circ -1}.
\end{align*}
Using Theorem~\ref{thm:mixprocomp_right_act_monoid},
\begin{align*}
\left(c\bar{@}d_1\right)\bar{@}d_2 &= c \mixprocomp \de{}\Bigg[\left(d_1^{-1} \mixcomp c\right) \shuffle
\Bigg(\left(d_{2}^{-1}\mixcomp c\right) \mixprocomp \\
& \quad \quad \left(\de{\left(d_1^{-1}\mixcomp c\right)^{\circ -1}} \circ \de{\left(d_1^{-1}\mixcomp c\right)}\right)\Bigg)\Bigg]^{\circ -1} \\
&= c \mixprocomp \de{\left(\left(d_{1}^{-1}\mixcomp c\right) \shuffle \left(\left(d_2^{-1} \mixcomp c\right) \mixprocomp \de{\doubleone}\right)\right)^{\circ -1}} \\
&= c \mixprocomp \de{\left(\left(d_{1}^{-1}\mixcomp c\right) \shuffle \left(d_{2}^{-1}\mixcomp c\right)\right)^{\circ -1}}.
\end{align*}
In light of Theorem~\ref{thm:Cauchy_prod_WF},
\begin{align*}
\left(c\bar{@}d_1\right)\bar{@}d_2 &= c \mixprocomp \de{\left(\left(d_1^{-1}.d_2^{-1}\right)\mixcomp c\right)^{\circ -1}} \\
&= c \mixprocomp \de{\left(\left(d_1.d_2\right)^{-1}\mixcomp c\right)^{\circ -1}} \\
&= c \bar{@}\left(d_1.d_2\right).
\end{align*}
\endpr 

It is important to note that for additive static output feedback product, known as {\em Wiener-Fliess feedback product}, the transformation group is the additive group $(\allcommutingseriesX{m},+,0)$ while for multiplicative static output feedback the multiplicative group $(\allcommutingpiseries{m},\cdot,\doubleone)$ performs the role.  

\section{Conclusions and Future work}
It was shown that the closed-loop system of a plant in Chen-Fliess series description in multiplicative output feedback with another system, given by Chen-Fliess series, has a Chen-Fliess series representation. An explicit expression of the closed-loop generating series was derived and the multiplicative dynamic feedback connection has a natural interpretation as a transformation group acting on the plant. It was then shown that when the Chen-Fliess series in the feedback is replaced by a memoryless map then the closed-loop system has a Chen-Fliess series representation. An explicit formula was provided for computing the generating series of the closed-loop system and it was shown that the multiplicative static feedback connection has a natural interpretation as a transformation group acting on the plant. Future work will be to address the solemn problem regarding the local convergence of the both multiplicative dynamic and static output feedback connections and to identify both the multiplicative dynamic and static feedback invariants.

\end{document}